# A GOLDEN PRODUCT IDENTITY FOR $e$

## ROBERT P. SCHNEIDER

ABSTRACT. We prove an infinite product representation for the constant $e$ involving the golden ratio, the Möbius function and the Euler phi function—prominent players in number theory.

Consider the number $e$, the base of the natural logarithm, a constant of fundamental significance in calculus and the sciences. A spectacular connection between $e$ and $\pi$ is arrived at through the theory of complex numbers, highlighted by Euler's identity [1]

$$e^{i\pi} = -1.$$

In this note, we relate $e$ to another famous constant via functions from number theory. Let $\tau$ denote the golden ratio (often denoted by $\phi$ in the popular literature), the number

$$\tau = \frac{1 + \sqrt{5}}{2}.$$

Discovered in antiquity, the golden ratio has fascinated mathematicians and laypersons alike over the centuries [3]. Among many nice properties that $\tau$ possesses, it is related to its own reciprocal by the identity $\tau = 1 + 1/\tau$.

In the theory of numbers, the Euler phi function $\varphi(n)$ counts the number of positive integers less than an integer $n$ that are co-prime to $n$. The classical Möbius function $\mu(n)$ is defined at a positive integer $n$ as

$$\mu(n) = \begin{cases} 1 \text{ if } n = 1 \\ 0 \text{ if } n \text{ is not squarefree} \\ 1 \text{ if } n \text{ is squarefree, having an even number of prime factors} \\ -1 \text{ if } n \text{ is squarefree, having an odd number of prime factors}. \end{cases}$$

Both $\varphi$ and $\mu$ are at the heart of deep results [4]. These two arithmetic functions, intimately connected to the factorization of integers, are related through the identity [2]

$$\varphi(n) = n \sum_{d|n} \frac{\mu(d)}{d}.$$

Here we prove an interesting infinite product expansion relating the phi function and Möbius function to the golden ratio; moreover, we show that the constant $e$ can be written in terms of these three number-theoretic objects. In the course of the proof, we find a curious inverse relationship between $\varphi$ and $\mu$.

**Theorem**. *We have the infinite product identity*

$$e = \prod_{n=1}^{\infty} \left(1 - \frac{1}{\tau^n}\right)^{\frac{\mu(n) - \varphi(n)}{n}}.$$

While the golden ratio is related to many marvelous self-similar phenomena, it is not particularly central to research in number theory. However, $e$, $\mu$, and $\varphi$ pervade the study of prime numbers and other arithmetic topics; it is intriguing to find these important entities intersecting $\tau$ in this product formula.

We begin the proof of the theorem by stating a few well-known identities. For any positive integer $n$, the following two identities hold [**2**]:

$$\sum_{d|n} \varphi(d) = n \qquad \sum_{d|n} \mu(d) = \begin{cases} 1 & \text{if } n = 1 \\ 0 & \text{if } n > 1. \end{cases}$$

In addition to these, we have the Maclaurin series

$$-\log(1-y) = \sum_{j=1}^{\infty} \frac{y^j}{j}$$

when $0 < y < 1$, where $\log x$ denotes the natural logarithm of a real variable $x > 0$. We use these expressions to prove two lemmas, leading to a quick proof of the main result.

**Lemma 1**. *On the interval* $0 < x < 1$, *we have the identities*

$$(1) \quad -\sum_{k=1}^{\infty} \frac{\varphi(k)}{k} \log(1-x^k) = \frac{x}{1-x}$$

$$(2) \quad -\sum_{k=1}^{\infty} \frac{\mu(k)}{k} \log(1-x^k) = x.$$

*Proof of Lemma 1*. To prove (1), observe that the left-hand side of the identity is equal to

$$\sum_{k=1}^{\infty} \frac{\varphi(k)}{k} \sum_{j=1}^{\infty} \frac{x^{kj}}{j}.$$

For any $n \geq 1$, the coefficient of $x^n$ in the preceding expression is

$$\frac{\sum_{k|n} \varphi(k)}{n} = 1.$$

Therefore the left side of (1) is equal to the geometric series

$$x + x^2 + x^3 + \cdots = \frac{x}{1-x}.$$

To prove (2), observe that the left-hand side of the identity is equal to

$$\sum_{k=1}^{\infty} \frac{\mu(k)}{k} \sum_{j=1}^{\infty} \frac{x^{kj}}{j}.$$

For any $n \geq 1$, the coefficient of $x^n$ in the preceding expression is

$$\frac{\sum_{k|n} \mu(k)}{n} = \begin{cases} 1 & \text{if } n = 1 \\ 0 & \text{if } n > 1. \end{cases}$$

Every coefficient but the first equals zero; therefore the left side of (2) is equal to $x$.

(These infinite series are absolutely convergent when $0 < x < 1$, justifying the changes in order of summation.) ∎

The proof of Lemma 1 makes use of basic properties of $\varphi$, $\mu$, geometric series, and power series. Using these identities, the following lemma displays a kind of inverse relationship between the Euler phi function and Möbius function with respect to the golden ratio: The two summations in Lemma 1 are reciprocals precisely when $x$ is equal to $1/\tau$.

**Lemma 2**. *We have the pair of reciprocal summation identities*

$$\tau = -\sum_{k=1}^{\infty} \frac{\varphi(k)}{k} \log\left(1 - \frac{1}{\tau^k}\right)$$

$$\frac{1}{\tau} = -\sum_{k=1}^{\infty} \frac{\mu(k)}{k} \log\left(1 - \frac{1}{\tau^k}\right).$$

*Proof of Lemma 2*. Observing that $0 < 1/\tau < 1$, we make the substitution $x = 1/\tau$ in the two identities of Lemma 1. We use the identity $1/\tau = \tau - 1$ to simplify the right-hand side of (1), as

$$\frac{1/\tau}{1 - 1/\tau} = \frac{1}{\tau - 1} = \frac{1}{1/\tau} = \tau. \quad \blacksquare$$

*Proof of the Theorem*. We subtract the second identity in Lemma 2 from the first identity and observe that $\tau - 1/\tau = 1$, to find

$$\sum_{k=1}^{\infty} \frac{\mu(k) - \varphi(k)}{k} \log\left(1 - \frac{1}{\tau^k}\right) = 1.$$

By applying the exponential function to both sides of this identity, we arrive at the theorem. ∎

Similar product identities for $\exp(x)$, $\exp\left(\frac{x}{1-x}\right)$, and $\exp\left(\frac{x^2}{1-x}\right)$ are obtained by applying the exponential function to the identities in Lemma 1. More generally, it is evident from the above proof that for any arithmetic function $f$, we have the identity

$$\prod_{n=1}^{\infty} (1 - x^n)^{-\frac{f(n)}{n}} = \exp\left(\sum_{n=1}^{\infty} \frac{(1 * f)(n)}{n} x^n\right),$$

where $(1 * f)(n) = \sum_{d|n} f(d)$ as usual.

## ACKNOWLEDGMENTS

The author thanks David Leep at the University of Kentucky for his guidance during the preparation of this report, and for noting a simplified proof of Lemma 1, as well as Cyrus Hettle, Richard Ehrenborg, Ken Ono, Adele Lopez, and Andrew Granville for their useful suggestions.

## REFERENCES

1. W. Dunham, *Euler: The Master of Us All*, The Mathematical Association of America, Washington, DC, 1999.
2. G. H. Hardy and E. M. Wright, *An Introduction to the Theory of Numbers*, Oxford University Press, Oxford, 1980.
3. M. Livio, *The Golden Ratio*, Broadway Books, New York, 2002.
4. J. Stopple, *A Primer of Analytic Number Theory: From Pythagoras to Riemann*, Cambridge University Press, Cambridge, 2003.

DEPARTMENT OF MATHEMATICS AND COMPUTER SCIENCE, EMORY UNIVERSITY, ATLANTA, GEORGIA 30322
*E-mail address:* `robert.schneider@emory.edu`